\documentclass[11pt,titlepage]{article}
\usepackage{amsfonts,amssymb,euscript,amscd,amsthm}
\usepackage[centertags]{amsmath}

\newtheorem {theorem} {Theorem}
\newtheorem {lemma} {Lemma}

\def \proof {\textbf {Proof. }}

\def \UU {\mathcal{U}}
\def \EEE {\mathcal{E}}
\def \R {\bf R}
\def \prt {\partial}
\def \ve {\varepsilon}
\def \E {\mbox{E}}
\def \suml {\sum \limits}
\def \intl {\int \limits}
\input{tcilatex}

\begin{document}

\title{Nonparametric Estimation in the Model of Moving Average. }
\author{Alexander Alekseev}
\maketitle



\vskip 2cm

\section{Introduction.}

The subject of robust estimation in time series is widely discussed in
literature. One of the approaches is to use GM-estimation. This method
incorporates a broad class of nonparametric estimators which under suitable
conditions includes estimators robust to outliers in data. For the linear
models the sensitivity of GM-estimators to outliers have been studied in the
work by Martin and Yohai [5], and influence functionals for this estimator
were derived. In this paper we follow this direction and examine the
asymptotical properties of the class of M-estimators, which is narrower than
the class of GM-estimators, but gives more insight into asymptotical
properties of such estimators. This paper gives an asymptotic expansion of
the residual weighted empirical process, which allows to prove asymptotic
normality of these estimators in case of non-smooth objective functions. For
simplicity MA(1) model is considered, but it will be shown that even in this
case mathematical techniques used to derive these asymptotic properties
appear to be rather complicated.However, the approach used in this paper
could be applied to GM-estimators and to more realistic models.

\section{Main Results.}

In this work we consider the model of moving average MA(1): 
\begin{equation*}
u_i = \ve_i - \alpha \ve_{i-1}, \quad i = 0,\pm 1, \pm 2, \dots, \leqno (1)
\end{equation*}
where $\{ \ve_i \}$ - iid, $E \ve_1=0,~E \ve_1^2<\infty$, $| \alpha | < 1$.%
\newline
Let $u_1, \dots, u_n$ be the observations of a random variable $u$.\newline
For every $\theta \in \R^1$ set 
\begin{equation*}
\ve_0 (\theta) = 0,
\end{equation*}
\begin{equation*}
\ve_i (\theta) = u_i + \theta \ve_{i-1} (\theta), \quad i = 1, 2, \dots %
\leqno (2)
\end{equation*}
It can easily be seen that 
\begin{equation*}
\ve_i (\theta) = \suml_{j=0}^{i-1} \theta^j u_{i-j}, \leqno (3)
\end{equation*}
where $\ve_i (\theta)$ are the residuals of the model (1). Assume for a
moment that the equation of a moving average holds only for $i = 1, 2, \dots$
with $\ve_0 = 0$.Also let's assume that variables $\{ \ve_i, \: i \ge 1 \}$
have a density function $g(x)$. Then the maximum likelihood equation for the
estimation of $\alpha$ can be constructed. Denote $\UU = (u_1, \dots,
u_n)^{T}, \quad \EEE = (\ve_1, \dots, \ve_n)^{T}$, and the matrix $J
(\alpha) $ as 
\begin{equation*}
J (\alpha) = \left( 
\begin{array}{ccccc}
1 & 0 & 0 & \ldots & 0 \\ 
\alpha & 1 & 0 & \ldots & 0 \\ 
\alpha^2 & \alpha & 1 & \ldots & 0 \\ 
\vdots & \vdots & \vdots & \ddots & \vdots \\ 
\alpha^{n-1} & \alpha^{n-2} & \alpha^{n-3} & \ldots & 1%
\end{array}
\right).
\end{equation*}
Then the equations $\displaystyle \ve_i = \suml_{j=0}^{i-1} \alpha^j
u_{i-j}, \quad i = 1, \dots, n$ can be rewritten as: 
\begin{equation*}
\EEE = J (\alpha) \UU
\end{equation*}
If $f_{\UU} (y_1, \dots, y_n)$ is the density function of the vector $\UU$,
then 
\begin{equation*}
f_{\UU} (y_1, \dots, y_n) = \prod \limits_{i=1}^n g \left( \suml_{j=0}^{i-1}
\alpha^j y_{i-j} \right).
\end{equation*}
Therefore, the maximum likelihood estimator, which is defined as a solution
of the maximization problem 
\begin{equation*}
\log f_{\UU} (y_1, \dots, y_n) \mathop {\rm \longrightarrow} \limits \sup
\limits_{\theta},
\end{equation*}
can be obtained as a root of the equation 
\begin{equation*}
\suml_{i=1}^n \frac {\prt \ve_i (\theta)} {\prt \theta} \frac {g^{\prime}(\ve%
_i (\theta))} {g (\ve_i (\theta))} = 0,
\end{equation*}
In this paper we examine a natural generalization of this estimator for the
model(1), namely, $M$-estimator $\widehat \alpha_n$ of the parameter $\alpha$%
. \newline
The estimator $\widehat \alpha_n$is defined as a solution of the equation 
\begin{equation*}
\frac {1} {n} \suml_{i=1}^n \frac {\prt \ve_i (\theta)} {\prt \theta} \Psi (%
\ve_i (\theta)) = 0,
\end{equation*}
with $\ve_i (\theta)$, determined by (2) or (3), where $\Psi (x)$ is an a
priori known function, which we will choose later.\newline
The asymptotic distribution of the estimator $\widehat \alpha_n$ in the case
when the function $\Psi (\cdot)$ satisfies $\mathop {\rm Var} \nolimits
|_{-\infty}^{+\infty} [\Psi] <\infty$ will be derived. This result will be
obtained with the help of asymptotically uniform expansion of the residual
weighted empirical process, which will be defined later and which will be of
interest by its own. \vskip 1cm

We fix some notation: 
\begin{equation*}
l_n (\theta) := \frac {1} {n} \suml_{k=1}^n \frac {\prt \ve_k (\theta)} {%
\prt \theta} \Psi (\ve_k (\theta)), \quad \widetilde l_n (\alpha) := \frac {1%
} {n} \suml_{k=1}^n \frac {\prt \ve_k (\alpha)} {\prt \theta} \Psi (\ve_k),
\end{equation*}
\begin{equation*}
u_n (x,\theta) := \frac {1} {n} \suml_{k=1}^n \frac {\prt \ve_k (\theta)} {%
\prt \theta} \mathop {\rm I} \nolimits (\ve_k (\theta) \le x), \quad
\widetilde u_n (x,\alpha) := \frac {1} {n} \suml_{k=1}^n \frac {\prt \ve_k
(\alpha)} {\prt \theta} \mathop {\rm I} \nolimits (\ve_k \le x),
\end{equation*}
where $\ve_i$ are derived from (1): 
\begin{equation*}
\ve_i = \suml_{j \ge 0} \alpha^j u_{i-j}.
\end{equation*}
If the variation of $\Psi$ is bounded,the following is true: 
\begin{equation*}
n^{1/2} \left[ l_n (\alpha + n^{-1/2} \tau) - \widetilde l_n (\alpha) \right]
= \intl_{-\infty}^{+\infty} \Psi (x) d n^{1/2} [ u_n (x,\alpha + n^{-1/2}
\tau) - \widetilde u_n (x,\alpha) ] =
\end{equation*}
\begin{equation*}
= -\intl_{-\infty}^{+\infty} n^{1/2} [ u_n (x,\alpha + n^{-1/2} \tau) -
\widetilde u_n (x,\alpha) ] d \Psi (x) + \left\{ n^{1/2} [ u_n (x,\alpha +
n^{-1/2} \tau) - \widetilde u_n (x,\alpha) ] \Psi (x)
\right\}_{-\infty}^{+\infty}.
\end{equation*}
It will be proved below that the second term is equal to zero. Now consider
the first term.

\begin{theorem}
Assume that the following conditions hold: 
\begin{equation*}
\E (\ve_1)^8 < \infty;
\end{equation*}
\begin{equation*}
g (x) > 0, \quad \lim \limits_{x \mathop {\rm \rightarrow} \limits \infty}
g(x) = 0, \quad \sup \limits_{x \in \R} | g^{\prime}(x) | < \infty.
\end{equation*}
Then 
\begin{equation*}
\sup \limits_{x \in \R, |\tau| < \theta} \left| n^{1/2} [ u_n (x,\alpha +
n^{-1/2} \tau) - \widetilde u_n (x,\alpha) ] + \tau g(x) \frac {\E \ve_1^2} {%
1-\alpha^2} \right| = o_{\mathop {\rm P} \nolimits} (1).
\end{equation*}
\end{theorem}

The first theorem will be proved in section 2. The next theorem is the main
result of the work and its proof utilizes the first theorem. \vskip 1.0cm

\begin{theorem}
Assume that the following conditions hold: 
\begin{equation*}
\E (\ve_1)^8 < \infty, \leqno (i)
\end{equation*}
\begin{equation*}
g (x) > 0, \quad \lim \limits_{x \mathop {\rm \rightarrow} \limits \infty} g
(x) = 0, \quad \lim \limits_{x \mathop {\rm \rightarrow} \limits -\infty} g
(x) = 0, \quad \sup \limits_{x \in \R} |g^{\prime}(x)| < \infty;
\end{equation*}
\begin{equation*}
\mathop {\rm Var} \nolimits |_{-\infty}^{+\infty} \, [\Psi] < \infty, \leqno %
(ii)
\end{equation*}
\begin{equation*}
\intl_{-\infty}^{\infty} g d \Psi \not = 0,
\end{equation*}
\begin{equation*}
E [\Psi(\ve_1)]=0.
\end{equation*}

then\newline
1) if $\Psi$ is continuous on $\R^1$, then with probability tending to one,
there exists a $n^{1/2}$--consistent solution of the following equation 
\begin{equation*}
n^{-1} \suml_{k=1}^n \frac {\prt \ve_k (\theta)} {\prt \theta} \Psi (\ve_k
(\theta)) = 0; \leqno (4)
\end{equation*}
2) for any $n^{1/2}$ - consistent solution $\widehat \alpha_n$ of(4) the
following statement

holds

\begin{equation*}
n^{1/2} (\widehat \alpha_n - \alpha) \mathop {\rm \longrightarrow}
\limits^{d} \mathop {\rm N} \nolimits (0, \sigma_{\Psi}^2 (\alpha)),
\end{equation*}
where 
\begin{equation*}
\sigma_{\Psi}^2 (\alpha) = (1-\alpha^2) \frac {\E \Psi^2 (\ve_1)} {\left( %
\intl_{-\infty}^{+\infty} g d \Psi \right)^2 \E (\ve_1)^2}.
\end{equation*}
\end{theorem}

\vskip 1.0cm

The function $\Psi (x) = F(x) - \frac {1} {2}$ where $F(x)$ is a continuous
distribution function of a certain zero mean symmetrical distribution
satisfies the conditions of theorem 2 in case $\ve_i$ are symmetrical zero
mean random variables. \vskip 1.5cm

\section{Proofs of the theorems}

\begin{lemma}
If for some $p \ge 1 \quad \E |\ve_1|^p < \infty$, then: 
\begin{equation*}
\sup \limits_{1 \le k \le n} \E \left| \frac {\prt \ve_k (\alpha)} {\prt %
\theta} \right|^p < \infty;
\end{equation*}
\begin{equation*}
\sup \limits_{1 \le k \le n} \E \left| \frac {\prt^2 \ve_k (\alpha)} {\prt %
\theta^2} \right|^p < \infty.
\end{equation*}
\end{lemma}

\proof 
\begin{equation*}
\frac {\prt \ve_k (\alpha)} {\prt \theta} = \suml_{j=1}^{k-1} j \alpha^{j-1}
u_{k-1}.
\end{equation*}
\begin{equation*}
\E \left| \frac {\prt \ve_k (\alpha)} {\prt \theta} \right|^p \le \E \left( %
\suml_{j=1}^{k-1} j |\alpha|^{j-1} |u_{k-1}| \right)^p \le c \E \left( \suml%
_{j=1}^{k-1} \alpha_1^j |u_{k-1}| \right)^p \le
\end{equation*}
\begin{equation*}
\le c \E \left( \suml_{j=1}^{k-1} \alpha_1^j |\ve_{k-1}| + |\alpha| \suml%
_{j=1}^{k-1} \alpha_1^j |\ve_{k-1-j}| \right)^p \le
\end{equation*}
\begin{equation*}
\le 2 c \left[ \E \left( \suml_{j=1}^{k-1} \alpha_1^j |\ve_{k-1}| \right)^p
+ |\alpha| \E \left( \suml_{j=1}^{k-1} \alpha_1^j |\ve_{k-1-j}| \right)^p %
\right],
\end{equation*}
where $\alpha_1 \in [0,+\infty): \quad \left\{ 
\begin{array}{c}
|\alpha| < \alpha_1 < 1 \\ 
j |\alpha|^{j-1} < c \alpha_1^j%
\end{array}
\right.$. Therefore, it is sufficient to prove 
\begin{equation*}
\sup \limits_{1 \le k \le n} \left\{ \E \left( \suml_{j=1}^{k-1} \alpha_1^j |%
\ve_{k-j}| \right)^p \right\}^{1/p} < \infty.
\end{equation*}
This clear from the Minkovsky's inequality:

\begin{equation*}
\left\{ \E \left( \suml_{j=1}^{k-1} \alpha_1^j |\ve_{k-j}| \right)^p
\right\}^{1/p} \le \alpha_1 \left\{ \E (\ve_1)^p \right\}^{1/p} + \left\{ \E %
\left( \suml_{j=2}^{k-1} \alpha_1^j |\ve_{k-j}| \right)^p \right\}^{1/p} \le
\dots \le
\end{equation*}
\begin{equation*}
\le (\alpha_1 + \dots + \alpha_1^{k-1}) \{ \E |\ve_1|^p \}^{1/p} < \infty,
\end{equation*}
what proves the first claim of the lemma. The second one can be proved in
the similar way. \vskip 1cm

\begin{lemma}
set 
\begin{equation*}
\sigma_k(\tau) := - \tau n^{-1/2} \suml_{t \ge k} \alpha^t \ve_{k-1-t} +
\tau n^{-1/2} \suml_{t=0}^{k-1} \left[ (\alpha + \tau n^{-1/2})^t - \alpha^t %
\right] \ve_{k-1-t} - \left( \alpha_n + \tau n^{-1/2} \right)^k \ve_0
\end{equation*}
if $\E(\ve_1)^4 < \infty$, then there exists such a $\widehat \sigma_k$, that

\begin{equation*}
\sup \limits_{|\tau| \le \theta} \left| \sigma_k (\tau) \right| \le \widehat
\sigma_k,
\end{equation*}
\begin{equation*}
\sup \limits_{1 \le k \le n} \E (\widehat \sigma_k)^4 < \infty.
\end{equation*}
\end{lemma}

\proof Let $B \in (0,1) : \quad |\alpha| + \theta n^{-1/2} < B$ for $n >
n_0, 0< \theta <\infty$, then $\sigma_k(\tau) \le \widehat \sigma_k$, where 
\begin{equation*}
\widehat \sigma_k := \theta n^{-1/2} \suml_{t \ge k} B^t |\ve_{k-1-t}| +
\theta^2 n^{-1} \suml_{t \ge 1} t B^{t-1} |\ve_{k-1-t}| - B^k \ve_0.
\end{equation*}
\begin{equation*}
\E (\sigma_k)^4 = \E {\left( \theta n^{-1/2} \suml_{t \ge k} B^t
|\ve_{k-1-t}| + \theta^2 n^{-1} \suml_{t \ge 1} t B^{t-1} |\ve_{k-1-t}| -
B^k \ve_0 \right)}^4.
\end{equation*}
For $\forall a,b > 0, p \ge 1$ we have the following inequality $(a+b)^p \le
2^{p-1} (a^p + b^p)$. Hence, it follows that 
\begin{equation*}
\E (\sigma_k)^4 \le 2^3 \cdot \left[ \E {\left( \theta n^{-1/2} \suml_{t \ge
k} B^t |\ve_{k-1-t}| \right)}^4 + \E {\left( \theta^2 n^{-1} \suml_{t \ge 1}
t B^{t-1} |\ve_{k-1-t}| - B^k \ve_0 \right)}^4 \right].
\end{equation*}
By Minkovsky's inequality: 
\begin{equation*}
{\left\{ \E {\left( \theta n^{-1/2} \suml_{t \ge k} B^t |\ve_{k-1-t}|
\right)}^4 \right\}}^{1/4} \le {\left\{ \theta^4 n^{-2} B^{4k} {\E (\ve_1)}%
^4 \right\}}^{1/4} +
\end{equation*}
\begin{equation*}
+ {\left\{ \E {\left( \theta n^{-1/2} \suml_{t \ge k+1} B^t |\ve_{k-1-t}|
\right)}^4 \right\}}^{1/4} \le {\left\{ \theta^4 n^{-2} \left( B^{4k} +
B^{4(k+1)} + \dots \right) \E (\ve_1)^4 \right\}}^{1/4} < \infty,
\end{equation*}
hence, 
\begin{equation*}
\E {\left( \theta n^{-1/2} \suml_{t \ge k} B^t |\ve_{k-1-t}| \right)}^4 <
\infty.
\end{equation*}
\begin{equation*}
{\left\{ \E {\left( \theta^2 n^{-1} \suml_{t \ge 1} t B^{t-1}
|\ve_{k-1-t}|-B^k \ve_0 \right)}^4 \right\}}^{1/4} \le
\end{equation*}
\begin{equation*}
\le {\left\{ \theta^8 n^{-4} \left( 1 \cdot B^0 + 2 \cdot B^1 + \dots + t
B^{t-1} + \dots \right) \E (\ve_1)^4-B^k E(\ve_1)^4 \right\}}^{1/4}.
\end{equation*}
Since the series $\suml_{t \ge 1} t B^{t-1}$ converges, it follows that 
\begin{equation*}
\E {\left( \theta^2 n^{-1} \suml_{t \ge 1} t B^{t-1} |\ve_{k-1-t}|-B^k \ve_0
\right)}^4 < \infty,
\end{equation*}
what proves the lemma. \vskip 1.0cm

\begin{lemma}
Let the following conditions hold: 
\begin{equation*}
\E (\ve_1)^8 < \infty \leqno (i)
\end{equation*}
$(ii) \quad g(x)$ is the density function of $\ve_i$ that satisfies the
following conditions: 
\begin{equation*}
\lim \limits_{x \mathop {\rm \rightarrow} \limits \infty} g(x) = 0,
\end{equation*}
\begin{equation*}
\sup \limits_{x} |g^{\prime}(x)| < \infty.
\end{equation*}
For $z_{1n} (x,\tau)$ defined as 
\begin{equation*}
z_{1n} (x,\tau) := n^{-1/2} \suml_{k=1}^n \left( \frac {\prt \ve_k (\alpha +
n^{-1/2} \tau)} {\prt \theta} - \frac {\prt \ve_k (\alpha)} {\prt \theta}
\right) \mathop {\rm I} \nolimits \left( \ve_k (\alpha + n^{-1/2} \tau) \le
x \right).
\end{equation*}
the following holds: 
\begin{equation*}
|z_{1n} (x,\tau)| = o_{\mathop {\rm P} \nolimits} (1).
\end{equation*}
\end{lemma}

\proof 
\begin{equation*}
z_{1n} (x,\tau) = n^{-1/2} \suml_{k=1}^n \left( \frac {\prt \ve_k (\alpha +
n^{-1/2} \tau)} {\prt \theta} - \frac {\prt \ve_k (\alpha)} {\prt \theta}
\right) \mathop {\rm I} \nolimits \left( \ve_k (\alpha + n^{-1/2} \tau) \le
x \right) =
\end{equation*}
\begin{equation*}
= n^{-1/2} \suml_{k=1}^n \tau n^{-1/2} \frac {\prt^2 \ve_k (\alpha)} {\prt %
\theta^2} \mathop {\rm I} \nolimits \left( \ve_k (\alpha + n^{-1/2} \tau)
\le x \right) +
\end{equation*}
\begin{equation*}
+ n^{-1/2} \suml_{k=1}^n \tau^2 n^{-1} \frac {\prt^3 \ve_k (\xi)} {\prt %
\theta^3} \mathop {\rm I} \nolimits \left( \ve_k (\alpha + n^{-1/2} \tau)
\le x \right) =
\end{equation*}
\begin{equation*}
= \tau n^{-1} \suml_{k=1}^n \frac {\prt^2 \ve_k (\alpha)} {\prt \theta^2} %
\mathop {\rm I} \nolimits \left( \ve_k (\alpha + n^{-1/2} \tau) \le x
\right) + o_{\mathop {\rm P} \nolimits} (1),
\end{equation*}
where $\xi \in [\alpha, \alpha + n^{-1/2} \tau]$. Therefore 
\begin{equation*}
|z_{1n} (x,\tau)| \le \left| \tau n^{-1} \suml_{k=1}^n \frac {\prt^2 \ve_k
(\alpha)} {\prt \theta^2} \left[ \mathop {\rm I} \nolimits \left( \ve_k
(\alpha + n^{-1/2} \tau) \le x \right) - \mathop {\rm I} \nolimits \left( \ve%
_k \le x \right) \right] \right| +
\end{equation*}
\begin{equation*}
+ \left| \tau n^{-1} \suml_{k=1}^n \frac {\prt^2 \ve_k (\alpha)} {\prt %
\theta^2} \mathop {\rm I} \nolimits \left( \ve_k \le x \right) \right| + o_{%
\mathop {\rm P} \nolimits} (1).
\end{equation*}
To transform the right hand of this inequality, we introduce the random
process $L_n^{+} (x,\tau)$: 
\begin{equation*}
L_n^{+} (x,\tau) := \tau n^{-1} \suml_{k=1}^n a_k^{+} \left[ \mathop {\rm I}
\nolimits \left( \ve_k (\alpha + n^{-1/2} \tau) \le x \right) - \mathop {\rm
I} \nolimits \left( \ve_k \le x \right) \right],
\end{equation*}
where 
\begin{equation*}
a_k^{+} := \left\{ 
\begin{array}{cl}
\frac {\prt^2 \ve_k (\alpha)} {\prt \theta^2}, & \mbox{if} \quad \frac {\prt%
^2 \ve_k (\alpha)} {\prt \theta^2} > 0, \\ 
0, & \mbox{if} \quad \frac {\prt^2 \ve_k (\alpha)} {\prt \theta^2} \le 0.%
\end{array}
\right.
\end{equation*}
\begin{equation*}
L_n^{-} (x,\tau) := \tau n^{-1} \suml_{k=1}^n a_k^{-} \left[ \mathop {\rm I}
\nolimits \left( \ve_k (\alpha + n^{-1/2} \tau) \le x \right) - \mathop {\rm
I} \nolimits \left( \ve_k \le x \right) \right],
\end{equation*}
\begin{equation*}
a_k^{-} := \left\{ 
\begin{array}{cl}
-\frac {\prt^2 \ve_k (\alpha)} {\prt \theta^2}, & \mbox{if} \quad \frac {\prt%
^2 \ve_k (\alpha)} {\prt \theta^2} < 0, \\ 
0, & \mbox{if} \quad \frac {\prt^2 \ve_k (\alpha)} {\prt \theta^2} \ge 0.%
\end{array}
\right.
\end{equation*}
\begin{equation*}
\left| \tau n^{-1} \suml_{k=1}^n \frac {\prt^2 \ve_k (\alpha)} {\prt \theta^2%
} \left[ \mathop {\rm I} \nolimits \left( \ve_k (\alpha + n^{-1/2} \tau) \le
x \right) - \mathop {\rm I} \nolimits \left( \ve_k \le x \right) \right]
\right| \le |L_n^{+} (x,\tau)| + |L_n^{-} (x,\tau)|.
\end{equation*}
We first show that $L_n^{+} (x,\tau) = o_{\mathop {\rm P} \nolimits} (1)$.
The proof utilizes the expansion obtained in lemma 7.1 [1] 
\begin{equation*}
\ve_k (\alpha + \tau n^{-1/2}) = \ve_k + \tau n^{-1/2} \mu_{k-1} + \sigma_k
(\tau),
\end{equation*}
and lemma 2. 
\begin{equation*}
| L_n^{+} (x,\tau) | = \left| \tau n^{-1} \suml_{k=1}^n a_k^{+} (\alpha) %
\left[ \mathop {\rm I} \nolimits \left( \ve_k \le x - \tau n^{-1/2}
\mu_{k-1} - \sigma_k (\tau) \right) - \mathop {\rm I} \nolimits (\ve_k \le x
) \right] \right| \le
\end{equation*}
\begin{equation*}
\le \left| \tau n^{-1} \suml_{k=1}^n a_k^{+} (\alpha) \left[ \mathop {\rm I}
\nolimits \left( \ve_k \le x - \tau n^{-1/2} \mu_{k-1} + \widehat \sigma_k
\right) - \mathop {\rm I} \nolimits (\ve_k \le x ) \right] \right|.
\end{equation*}
Consider random processes $v_n (x,\tau), \widetilde v_n (x)$: 
\begin{equation*}
v_n (x,\tau) := n^{-1} \suml_{k=1}^n a_k^{+} \left[ \mathop {\rm I}
\nolimits \left( \ve_k \le x - \tau n^{-1/2} \mu_{k-1} + \widehat \sigma_k
\right) - G \left( x - \tau n^{-1/2} \mu_{k-1} + \widehat \sigma_k \right) %
\right],
\end{equation*}
\begin{equation*}
\widetilde v_n (x) := n^{-1} \suml_{k=1}^n a_k^{+} \left[ \mathop {\rm I}
\nolimits \left( \ve_k \le x \right) - G \left( x \right) \right].
\end{equation*}
then 
\begin{equation*}
\left| \tau n^{-1} \suml_{k=1}^n a_k^{+} (\alpha) \left[ \mathop {\rm I}
\nolimits \left( \ve_k \le x - \tau n^{-1/2} \mu_{k-1} + \widehat \sigma_k
\right) - \mathop {\rm I} \nolimits (\ve_k \le x ) \right] \right| \le
\end{equation*}
\begin{equation*}
\le \left| \tau ( v_n (x,\tau) - \widetilde v_n (x) ) \right| + \left|
n^{-1} \suml_{k=1}^n a_k^{+} (\alpha) \left( G \left( x - \tau n^{-1/2}
\mu_{k-1} + \widehat \sigma_k \right) - G ( x ) \right) \right|,
\end{equation*}
and 
\begin{equation*}
\left| n^{-1} \suml_{k=1}^n a_k^{+} (\alpha) \left( G \left( x - \tau
n^{-1/2} \mu_{k-1} + \widehat \sigma_k \right) - G ( x ) \right) \right| =
\end{equation*}
\begin{equation*}
= \left| n^{-3/2} \suml_{k=1}^n a_k^{+} (\alpha) \tau \mu_{k-1} g(\xi) +
n^{-1} \suml_{k=1}^n a_k^{+} (\alpha) \widehat \sigma_k g(\xi) + o_{\mathop
{\rm P} \nolimits} (1) \right| = o_{\mathop {\rm P} \nolimits} (1),
\end{equation*}
where $\xi \in [x - \tau n^{-1/2} \mu_{k-1} + \widehat \sigma_k, x]$. To
prove that $\left| v_n (x,\tau) - \widetilde v_n (x) \right| = o_{\mathop 
\mathrm{P} \nolimits} (1)$ we use theorem 2.1 from [2]. Let's check whether
the conditions of this theorem . The condition $\E | \ve_1 |^8 < \infty$
implies $\E | \ve_1 |^4 < \infty$, and because of the lemma1 the latter
yields 
\begin{equation*}
\sup \limits_{1 \le k \le n} \E \left( \frac {\prt^2 \ve_k (\alpha)} {\prt %
\theta^2} \right)^4 < \infty.
\end{equation*}
Hence, 
\begin{equation*}
n^{-1} \suml_{k=1}^n \E \left( a_k^{+} (\alpha) \right)^4 = O (1).
\end{equation*}
To check the condition $n^{-1} \suml_{k=1}^n \E \left[ \left( a_k^{+}
(\alpha) \right)^4 {\widehat \sigma_k}^2 \right] = O (1)$ we use
Cauchy-Bunyakovskii inequality: 
\begin{equation*}
n^{-1} \suml_{k=1}^n \E \left[ \left( a_k^{+} (\alpha) \right)^4 {\widehat
\sigma_k}^2 \right] \le n^{-1} \suml_{k=1}^n \left[ \E \left( a_k^{+}
(\alpha) \right)^8 \right]^{1/2} \left[ \E {\widehat \sigma_k}^4 \right]%
^{1/2}.
\end{equation*}
because of the lemma1 and lemma2: 
\begin{equation*}
\sup_{1 \le k \le n} \E \left( a_k^{+} (\alpha) \right)^8 < \infty,
\end{equation*}
\begin{equation*}
\sup_{1 \le k \le n} \E {\widehat \sigma_k}^4 < \infty,
\end{equation*}
therefore, 
\begin{equation*}
n^{-1} \suml_{k=1}^n \E \left[ \left( a_k^{+} (\alpha) \right)^4 {\widehat
\sigma_k}^2 \right] = O (1).
\end{equation*}
Now we check the following condition: 
\begin{equation*}
n^{-1} \suml_{k=1}^n a_k^{+} (\alpha) {\widehat \sigma_k}^2 = O_{\mathop
{\rm P} \nolimits} (1).
\end{equation*}
Lemma 2 yields $\sup \limits_{1 \le k \le n} \E {\widehat \sigma_k}^2 <
\infty$, 
\begin{equation*}
n^{-1} \suml_{k=1}^n a_k^{+} (\alpha) {\widehat \sigma_k}^2 \le \left( \sup
\limits_{1 \le k \le n} {\widehat \sigma_k}^2 \right) n^{-1} \suml_{k=1}^n
a_k^{+} (\alpha) = \left( \sup \limits_{1 \le k \le n} {\widehat \sigma_k}^2
\right) n^{-1} \suml_{k=1}^n \left| \frac {\prt^2 \ve_k (\alpha)} {\prt %
\theta^2} \right|.
\end{equation*}
Consider a random process 
\begin{equation*}
\frac {\prt^2 \widetilde \ve (\alpha)} {\prt \theta^2} = \suml%
_{j=2}^{\infty} \left( j (j-1) | \alpha |^{j-2} | u_{-j} | \right).
\end{equation*}
then 
\begin{equation*}
\E \left( \frac {\prt^2 \widetilde \ve (\alpha)} {\prt \theta^2} - \left| 
\frac {\prt^2 \ve_k (\alpha)} {\prt \theta^2} \right| \right)^2 \le \E %
\left( \suml_{j=k+1}^{\infty} j (j-1) | \alpha |^{j-2} | u_{-j} | \right)^2
\le \E \left( \suml_{j=k+1}^{\infty} c \alpha^j | u_{-j} | \right)^2 \le
\end{equation*}
\begin{equation*}
\le 2 c^2 \left[ \E \left( \suml_{j=k+1}^{\infty} \alpha_1^j | \ve_{-j} |
\right)^2 + \E \left( | \alpha| \suml_{j=k+1}^{\infty} \alpha_1^j | \ve%
_{-j-1} | \right)^2 \right] =
\end{equation*}
\begin{equation*}
= 2 c^2 \E \ve_1^2 \left[ \suml_{j=k+1}^{\infty} \alpha_1^{2j} + | \alpha | %
\suml_{j=k+1}^{\infty} \alpha_1^{2j} \right] \mathop {\rm \longrightarrow}
\limits 0, k \mathop {\rm \rightarrow} \limits \infty,
\end{equation*}
where $\alpha_1: \left\{ 
\begin{array}{l}
|\alpha| < \alpha_1 < 1 \\ 
j (j-1) |\alpha|^{j-2} \le c \alpha_1^j%
\end{array}
\right. $.\newline
Hence 
\begin{equation*}
n^{-1} \suml_{k=1}^n \left| \frac {\prt^2 \ve_k (\alpha)} {\prt \theta^2}
\right| \mathop {\rm \longrightarrow} \limits^{\mathop {\rm P} \nolimits} \E %
\left( \suml_{j=2}^{\infty} j (j-1) | \alpha |^{j-2} | u_{-j} | \right)^2,
\end{equation*}
and 
\begin{equation*}
n^{-1} \suml_{k=1}^n a_k^{+} (\alpha) {\widehat \sigma_k}^2 = O_{\mathop
{\rm P} \nolimits} (1).
\end{equation*}
So all the conditions of the theorem 2.1 from [2]hold, therefore

\begin{equation*}
\sup \limits_{x \in \R^1, |\tau| \le \theta} n^{1/2} \left[ v_n (x,\tau) -
\widetilde v_n (x) \right] = o_{\mathop {\rm P} \nolimits} (1).
\end{equation*}
we have proved that 
\begin{equation*}
\left| \tau n^{-1} \suml_{k=1}^n a_k^{+} (\alpha) \left[ \mathop {\rm I}
\nolimits \left( \ve_k \le x - \tau n^{-1/2} \mu_{k-1} + \widehat \sigma_k
\right) - \mathop {\rm I} \nolimits (\ve_k \le x ) \right] \right| = o_{%
\mathop {\rm P} \nolimits} (1),
\end{equation*}
and 
\begin{equation*}
| L_n^{+} (x,\tau) | = o_{\mathop {\rm P} \nolimits} (1).
\end{equation*}
Therefore, 
\begin{equation*}
\left| \tau n^{-1} \suml_{k=1}^n \frac {\prt^2 \ve_k (\alpha)} {\prt \theta^2%
} \left[ \mathop {\rm I} \nolimits \left( \ve_k (\alpha + n^{-1/2} \tau) \le
x \right) - \mathop {\rm I} \nolimits \left( \ve_k \le x \right) \right]
\right| = o_{\mathop {\rm P} \nolimits} (1).
\end{equation*}
To finish the proof of the lemma we need to show that 
\begin{equation*}
\left| \tau n^{-1} \suml_{k=1}^n \frac {\prt^2 \ve_k (\alpha)} {\prt \theta^2%
} \mathop {\rm I} \nolimits \left( \ve_k \le x \right) \right| = o_{\mathop
{\rm P} \nolimits} (1).
\end{equation*}
\begin{equation*}
\left| \tau n^{-1} \suml_{k=1}^n \frac {\prt^2 \ve_k (\alpha)} {\prt \theta^2%
} \mathop {\rm I} \nolimits \left( \ve_k \le x \right) \right| \le \left|
\tau n^{-1} \suml_{k=1}^n \frac {\prt^2 \ve_k (\alpha)} {\prt \theta^2} %
\left[ \mathop {\rm I} \nolimits \left( \ve_k \le x \right) - G (x) \right]
\right| + \left| \tau n^{-1} G (x) \suml_{k=1}^n \frac {\prt^2 \ve_k (\alpha)%
} {\prt \theta^2} \right|.
\end{equation*}
We apply theorem 2.1 from [2] to process $n^{-1} \suml_{k=1}^n \frac {\prt^2 %
\ve_k (\alpha)} {\prt \theta^2} \left[ \mathop {\rm I} \nolimits \left( \ve%
_k \le x \right) - G (x) \right]$\newline
In the same way it can be shown that 
\begin{equation*}
n^{-1} \suml_{k=1}^n \frac {\prt^2 \ve_k (\alpha)} {\prt \theta^2} \mathop
{\rm \longrightarrow} \limits^{\mathop {\rm P} \nolimits} \E \left( \suml%
_{j=2}^{\infty} j (j-1) \alpha^{j-2} u_{-j} \right) = 0.
\end{equation*}
This finishes the proof of the lemma. \vskip 1.0cm

\begin{lemma}
Assume that the following conditions hold: 
\begin{equation*}
\E (\ve_1)^8 < \infty; \leqno (i)
\end{equation*}
\begin{equation*}
g (x) > 0, \quad \lim \limits_{x \mathop {\rm \rightarrow} \limits \infty} g
(x) = 0, \quad \sup \limits_{x} | g^{\prime}(x) | < \infty \leqno (ii).
\end{equation*}
Let 
\begin{equation*}
z_{2n} (x,\tau) := n^{-1/2} \suml_{k=1}^n \frac {\prt \ve_k (\alpha)} {\prt %
\theta} \left[ \mathop {\rm I} \nolimits \left( \ve_k (\alpha + n^{-1/2}
\tau) \le x \right) - \mathop {\rm I} \nolimits \left( \ve_k \le x \right) %
\right].
\end{equation*}
then 
\begin{equation*}
\left| z_{2n} (x,\tau) + \tau g (x) \frac {1} {1 - \alpha^2} \E \ve_1^2
\right| = o_{\mathop {\rm P} \nolimits} (1).
\end{equation*}
\end{lemma}

\proof 
\begin{equation*}
z_{2n} (x,\tau) = n^{-1/2} \suml_{k=1}^n \frac {\prt \ve_k (\alpha)} {\prt %
\theta} \left[ \mathop {\rm I} \nolimits \left( \ve_k \le x - \tau n^{-1/2}
\mu_{k-1} - \sigma_k (\tau) \right) - \mathop {\rm I} \nolimits \left( \ve_k
\le x \right) \right].
\end{equation*}
consider a process $v_n (x,\tau)$: 
\begin{equation*}
v_n (x,\tau) := n^{-1} \suml_{k=1}^n \frac {\prt \ve_k (\alpha)} {\prt \theta%
} \left\{ \mathop {\rm I} \nolimits \left( \ve_k \le x - \tau n^{-1/2}
\mu_{k-1} - \sigma_k (\tau) \right) - G (x - \tau n^{-1/2} \mu_{k-1} -
\sigma_k (\tau)) \right\},
\end{equation*}
then 
\begin{equation*}
z_{2n} = n^{1/2} \left[ v_n (x,\tau) - \widetilde v_n (x) \right] + n^{-1/2} %
\suml_{k=1}^n \frac {\prt \ve_k (\alpha)} {\prt \theta} \left[ G (x - \tau
n^{-1/2} \mu_{k-1} - \sigma_k (\tau)) - G (x) \right].
\end{equation*}
By analogy to lemma 3 it can be shown that 
\begin{equation*}
\sup \limits_{x \in {\mathbf{R}}^1, |\tau| \le \theta} n^{1/2} \left[ v_n
(x,\tau) - \widetilde v_n (x) \right] = o_{\mathop {\rm P} \nolimits} (1).
\end{equation*}
Consider a process 
\begin{equation*}
n^{-1/2} \suml_{k=1}^n \frac {\prt \ve_k (\alpha)} {\prt \theta} \left[ G (x
- \tau n^{-1/2} \mu_{k-1} - \sigma_k (\tau)) - G (x) \right] =
\end{equation*}
\begin{equation*}
= n^{-1/2} \suml_{k=1}^n \frac {\prt \ve_k (\alpha)} {\prt \theta} \left( -
\tau n^{-1/2} \mu_{k-1} - \sigma_k (\tau) \right) g (x) +
\end{equation*}
\begin{equation*}
+ n^{-1/2} \suml_{k=1}^n \frac {\prt \ve_k (\alpha)} {\prt \theta} \left( -
\tau n^{-1/2} \mu_{k-1} - \sigma_k (\tau) \right)^2 g^{\prime}(\xi) \mathop
{\rm \longrightarrow} \limits^{\mathop {\rm P} \nolimits}
\end{equation*}
\begin{equation*}
\mathop {\rm \longrightarrow} \limits^{\mathop {\rm P} \nolimits} -\tau g
(x) \E \left\{ \left( \suml_{j=1}^{\infty} j \alpha^{j-1} u_{k-j} \right)
\left( \suml_{j=0}^{\infty} \alpha^j \ve_{k-j-1} \right) \right\},
\end{equation*}
where $\xi \in [x - \tau n^{-1/2} \mu_{k-1} - \sigma_k (\tau), x]$, and 
\begin{equation*}
\E \left\{ \left( \suml_{j=1}^{\infty} j \alpha^{j-1} u_{k-j} \right) \left( %
\suml_{j=0}^{\infty} \alpha^j \ve_{k-j-1} \right) \right\} = \frac {1} {%
1-\alpha^2} \E \ve_1^2,
\end{equation*}
what proves the lemma.\newline
\vskip 1.0cm \textbf{Proof of Theorem 1}\newline
Rewrite the process in the following way: 
\begin{equation*}
n^{1/2} [ u_n (x,\alpha + n^{-1/2} \tau) - u_n (x,\alpha) ] = z_{1n}
(x,\tau) + z_{2n} (x,\tau),
\end{equation*}
where 
\begin{equation*}
z_{1n} (x,\tau) := n^{-1/2} \suml_{k=1}^n \left( \frac {\prt \ve_k (\alpha +
n^{-1/2} \tau)} {\prt \theta} - \frac {\prt \ve_k (\alpha)} {\prt \theta}
\right) \mathop {\rm I} \nolimits \left( \ve_k (\alpha + n^{-1/2} \tau) \le
x \right),
\end{equation*}
\begin{equation*}
z_{2n} (x,\tau) := n^{-1/2} \suml_{k=1}^n \frac {\prt \ve_k (\alpha)} {\prt %
\theta} \left[ \mathop {\rm I} \nolimits \left( \ve_k (\alpha + n^{-1/2}
\tau) \le x \right) - \mathop {\rm I} \nolimits \left( \ve_k \le x \right) %
\right].
\end{equation*}
Now it can be seen that lemma 3 and lemma4 prove the theorem.\newline
\vskip 1.0cm \textbf{Proof of Theorem 2}\newline
1) 
\begin{equation*}
\lambda (\alpha) = - \intl_{-\infty}^{\infty} g d \Psi \cdot \E \ve_1^2
\cdot \frac {1} {1-\alpha^2},
\end{equation*}
\begin{equation*}
l_n (\theta)= n^{-1} \suml_{k=1}^n \frac {\prt \ve_k (\theta)} {\prt \theta}
\Psi (\ve_k (\theta)), \quad \widetilde l_n (\alpha) := \frac {1} {n} \suml%
_{k=1}^n \frac {\prt \ve_k (\alpha)} {\prt \theta} \Psi (\ve_k)
\end{equation*}

Without loss of generality we consider $\lambda (\alpha) > 0$, such that
there exists $A > 0$ satisfying 
\begin{equation*}
n^{1/2} l_n (\alpha + n^{-1/2} A) = n^{1/2} \widetilde l_n (\alpha) +
\lambda (\alpha) A + o_{\mathop {\rm P} \nolimits} (1) > 0; \leqno (2)
\end{equation*}
\begin{equation*}
n^{1/2} l_n (\alpha - n^{-1/2} A) = n^{1/2} \widetilde l_n (\alpha) -
\lambda (\alpha) A + o_{\mathop {\rm P} \nolimits} (1) < 0. \leqno (3)
\end{equation*}
Since $l_n (\theta)$ is continuous on $\theta$ in some neighbourhood of $%
\alpha$, then because of (2) and (3), a $n^{1/2}$~- consistent solution of
(1)exists in interval $(\alpha - n^{-1/2} A; \alpha + n^{-1/2} A)$. \newline
2) 
\begin{equation*}
n^{1/2} \left[ l_n (\alpha + n^{-1/2} \tau) - \widetilde l_n (\alpha) \right]
= - \intl_{-\infty}^{\infty} n^{1/2} \left[ u_n (x, \alpha + n^{-1/2} \tau)
- \widetilde u_n (x, \alpha) \right] d \Psi (x) +
\end{equation*}
\begin{equation*}
+ \left\{ n^{1/2} \left[ u_n (x, \alpha + n^{-1/2} \tau) - \widetilde u_n
(x, \alpha) \right] \Psi (x) \right\}_{-\infty}^{+\infty}.
\end{equation*}
Since $\mathop {\rm Var} \nolimits |_{-\infty}^{+\infty} \, [\Psi (x)] <
\infty$ holds , and because of theorem 1 and conditions imposed on $g (x)$,
it can be obtained that $\left\{ n^{1/2} \left[ u_n (x, \alpha + n^{-1/2}
\tau) - \widetilde u_n (x, \alpha) \right] \Psi (x)
\right\}_{-\infty}^{+\infty} = 0$, and also that 
\begin{equation*}
n^{1/2} \left[ l_n (\alpha + n^{-1/2} \tau) - \widetilde l_n (\alpha) \right]
= - \intl_{-\infty}^{\infty} n^{1/2} \left[ u_n (x, \alpha + n^{-1/2} \tau)
- \widetilde u_n (x, \alpha) \right] d \Psi (x) =
\end{equation*}
\begin{equation*}
= \tau \cdot \frac {1} {1-\alpha^2} \cdot \E \ve_1^2 \cdot \intl%
_{-\infty}^{+\infty} g d \Psi + o_{\mathop {\rm P} \nolimits} (1)
\end{equation*}
uniformly on $|\tau| \le \Theta$.\newline
Let $\widehat \alpha_n$ be a $n^{1/2}$-consistent solution of (1), that is $%
n^{1/2} (\widehat \alpha_n - \alpha) = O_{\mathop {\rm P} \nolimits} (1)$.
Since 
\begin{equation*}
n^{1/2} l_n (\widehat \alpha_n) = n^{1/2} \widetilde l_n (\alpha) + \lambda
(\alpha) n^{1/2} (\widehat \alpha_n - \alpha) + o_{\mathop {\rm P}
\nolimits} (1) = 0,
\end{equation*}
then 
\begin{equation*}
n^{1/2} (\widehat \alpha_n - \alpha) = - [\lambda (\alpha)]^{-1} \cdot
n^{1/2} \widetilde l_n (\alpha) + o_{\mathop {\rm P} \nolimits} (1).
\end{equation*}
Consider the process 
\begin{equation*}
n^{1/2} \widetilde l_n (\alpha) = n^{-1/2} \suml_{k=1}^n \frac {\prt \ve_k
(\theta)} {\prt \theta} \cdot \Psi (\ve_k).
\end{equation*}
If we denote 
\begin{equation*}
\frac {\prt \widetilde \ve_k (\alpha)} {\prt \theta} = \suml_{j=1}^{\infty}
j \alpha^{j-1} u_{k-j},
\end{equation*}
then 
\begin{equation*}
\E \left( \frac {\prt \widetilde \ve_k (\alpha)} {\prt \theta} - \frac {\prt %
\ve_k (\alpha)} {\prt \theta} \right)^2 = \E \left( \suml_{j=k+1}^{\infty} j
\alpha^{j-1} u_{k-j}\right)^2 \le c \alpha_1^k
\end{equation*}
for some $\alpha_1, \quad 1 > \alpha_1 > |\alpha|$. 
\begin{equation*}
\E \left| n^{-1/2} \suml_{k=1}^n \left( \frac {\prt \widetilde \ve_k (\alpha)%
} {\prt \theta} - \frac {\prt \ve_k (\alpha)} {\prt \theta} \right) \Psi (\ve%
_k) \right| \le
\end{equation*}
\begin{equation*}
\le n^{-1/2} \suml_{k=1}^n \left\{ \E \left( \frac {\prt \widetilde \ve_k
(\alpha)} {\prt \theta} - \frac {\prt \ve_k (\alpha)} {\prt \theta}
\right)^2 \right\}^{1/2} \cdot \left\{ \E \Psi^2 (\ve_k) \right\}^{1/2} \le
\end{equation*}
\begin{equation*}
\le c n^{-1/2} \suml_{k=1}^n \alpha_1^k \left\{ \E \Psi^2 (\ve_k)
\right\}^{1/2} = c n^{-1/2} \left( \E \Psi^2 (\ve_1) \right) \suml_{k=1}^n
\alpha_1^k = o (1). \eqno (*)
\end{equation*}
Since $\displaystyle \xi_k := \frac {\prt \widetilde \ve_k (\alpha)} {\prt %
\theta} \Psi (\ve_k)$ is strictly stationary and forms a
martingale-difference, then the central limit theorem can be applied to it:

\begin{equation*}
n^{-1/2} \suml_{k=1}^n \xi_k \mathop {\rm \longrightarrow} \limits \mathop
{\rm N} \nolimits \left( 0, \E \left( \suml_{j=1}^{\infty} j \alpha^{j-1}
u_{-j} \right)^2 \cdot \E \Psi^2 (\ve_1) \right).
\end{equation*}
Because of the proved property (*) we obtain

\begin{equation*}
n^{1/2} \widetilde l_n (\alpha) \mathop {\rm \longrightarrow} \limits %
\mathop {\rm N} \nolimits \left( 0, \E \left( \suml_{j=1}^{\infty} j
\alpha^{j-1} u_{-j} \right)^2 \cdot \E \Psi^2 (\ve_1) \right).
\end{equation*}
\begin{equation*}
n^{1/2} (\widehat \alpha_n - \alpha) = - [\lambda (\alpha)]^{-1} \cdot
n^{1/2} \widetilde l_n (\alpha) + o_{\mathop {\rm P} \nolimits} (1) \mathop
{\rm \longrightarrow} \limits \mathop {\rm N} \nolimits \left( 0, \frac {\E %
\left( \suml_{j=1}^{\infty} j \alpha^{j-1} u_{-j} \right)^2 \cdot \E \Psi^2 (%
\ve_1)} {\left( \intl_{-\infty}^{+\infty} g d \Psi \right)^2 \cdot \frac {1%
} {(1-\alpha^2)^2} (\E \ve_1^2)^2} \right).
\end{equation*}
Since 
\begin{equation*}
\E \left( \suml_{j=1}^{\infty} j \alpha^{j-1} u_{-j} \right)^2 = \frac {\E %
\ve_1^2} {1-\alpha^2},
\end{equation*}
then 
\begin{equation*}
\sigma_{\Psi}^2 (\alpha) = (1-\alpha^2) \frac {\E \Psi^2 (\ve_1)} {\E \ve%
_1^2 \left( \intl_{-\infty}^{+\infty} g d \Psi \right)^2},
\end{equation*}
what proves theorem 2. 

\end{document}